\documentclass[12pt]{article}
\usepackage{amsmath,amssymb,fullpage}
\usepackage[applemac]{inputenc}

\newenvironment{myenumerate}{%

\begin{enumerate}}{\end{enumerate}}

\newcommand{\dproof}{\noindent {Proof.} \quad}
\newcommand{\fproof}{\hfill $\square$ \bigskip}

\newtheorem{definition}{Definition}[section]

\newtheorem{theorem}[definition]{Theorem}
\newtheorem{problem}[definition]{Problem}

\numberwithin{equation}{section}

\def\RB{\mathbb{R}}

\def\FB{\mathbb{F}}
\def\EB{\mathbb{E}}
\def\UB{\mathbb{U}}

\def\GB{\mathbb{G}}

\def\GC{\mathcal{G}}
\def\FC{\mathcal{F}}

\def\AC{\mathcal{A}}

\def\RC{\mathcal R}

\def\ess{\mathop{ess \;sup}}

\def\1B{\text{1\!\!I}}

\def\tY{\tilde{Y}}
\def\tN{\tilde{N}}

\def\tb{\tilde{b}}

\def\tZ{\tilde{Z}}
\def\tK{\tilde{K}}

\def\tX{\tilde{X}}
\def\tsi{\tilde{\sigma}}
\def\tga{\tilde{\gamma}}
\def\tg{\tilde{g}}
\def\tA{\tilde{A}}

\def\hX{\hat{X}}
\def\hZ{\hat{Z}}
\def\hA{\hat{A}}
\def\hY{\hat{Y}}

\def\hu{\hat{u}}

\def\hb{\hat{b}}
\def\hf{\hat{f}}

\def\hK{\hat{K}}
\def\hH{\hat{H}}
\def\hp{\hat{p}}

\def\hq{\hat{q}}
\def\hr{\hat{r}}

\def\hla{\hat{\lambda}}

\begin{document}

\title{Optimal control of predictive mean-field equations and applications to finance}

\date{19 May 2015}
\author{
Bernt \O ksendal$^{1,2,3}$ \and Agn\`es Sulem$^{4,5}$}

\footnotetext[1]{
Dept. of Mathematics, University of Oslo,
P.O. Box 1053 Blindern, N--0316 Oslo, Norway,
email: {\tt oksendal@math.uio.no}.}
\footnotetext[2]{
Norwegian School of Economics, Helleveien 30, N--5045 Bergen, Norway.
}
\footnotetext[3]{This research was carried out with support of CAS - Centre for Advanced Study, at the Norwegian Academy of Science and Letters, within the research program SEFE.}
\footnotetext[4]{ INRIA Paris-Rocquencourt, Domaine de Voluceau, Rocquencourt, BP 105, Le Chesnay Cedex, 78153, France,  and Universit\'e Paris-Est, email: {\tt agnes.sulem@inria.fr}}
\footnotetext[5]{ Universit\'e Paris-Est - Marne la Vall\'ee}

\maketitle



\begin{abstract}
We study a coupled system of controlled stochastic differential equations (SDEs) driven by a Brownian motion and a compensated Poisson random measure, consisting of a forward SDE in the unknown process $X(t)$ and a \emph{predictive mean-field} backward SDE (BSDE) in the unknowns $Y(t), Z(t), K(t,\cdot)$. The driver of the BSDE at time $t$ may depend not just upon the unknown processes $Y(t), Z(t), K(t,\cdot)$, but also on the predicted future value $Y(t+\delta)$, defined by the conditional expectation $A(t):= E[Y(t+\delta) | \mathcal{F}_t]$. \\
We give a sufficient and a necessary maximum principle for the optimal control of such systems, and then we apply these results to the following two problems:\\
(i)
Optimal portfolio in a financial market with an \emph{insider influenced asset price process.} \\
(ii) 
 Optimal consumption rate from a cash flow modeled as a geometric It\^ o-L\'  evy SDE, with respect to \emph{predictive recursive utility}.
\end{abstract}

\section{Introduction}\label{sec1}
The purpose of this paper is to introduce and study a pricing model where beliefs about the future development of the price process influence its current dynamics. We think this can be a realistic assumption in price dynamics where human psychology is involved, for example in electricity prices, oil prices and energy markets in general. It can also be a natural model of the risky asset price in an insider influenced market. See Section 5.1.\\
We model such price processes as backward stochastic differential equations (BSDEs) driven by Brownian motion and a compensated Poisson random measure, where the coefficients depend not only of the current values of the unknown processes, but also on their predicted future values. These predicted values are expressed mathematically in terms of conditional expectation, and we therefore name such equations \emph{predictive mean-field equations}. To the best of our knowledge such systems have never been studied before.

In applications to portfolio optimization in a financial market where the price process is modeled by a predictive mean-field equation, we are led to consider coupled systems of forward-backward stochastic differential equations (FBSEDs), where the BSDE is of predictive mean-field type. In this paper we study solution methods for the optimal control of such systems in terms of maximum principles. Then we apply these methods to study\\
(i)
optimal portfolio in a financial market with an \emph{insider influenced asset price process.} (Section 5.1), and  \\
(ii) 
 optimal consumption rate from a cash flow modeled as a geometric It\^ o-L\'  evy SDE, with respect to \emph{predictive recursive utility} (Section 5.2).

\section{Formulation of the problem}\label{sec2}
We now present our model in detail.

Let $B(t) = B(t,\omega); \; (t,\omega) \in [0, \infty) \times \Omega$ and $\tN(dt, d\zeta) = N(dt, d \zeta) - \nu (d\zeta)dt$ be a Brownian motion and an independent compensated Poisson random measure, respectively, on a filtered probability space  $\left(\Omega, \EB, \FB = \{ \FC_t\}_{t \geq 0}, P\right)$ satisfying the usual conditions. We consider a controlled system of  \emph{predictive} (time-advanced) coupled mean-field forward-backward stochastic differential equations (FBSDEs) of the form ($T > 0$  and $\delta>0$ are given constants)
\begin{itemize}
\item {Forward SDE in $X(t$):}
\begin{equation}\label{eq2.1a}
\begin{cases}
dX(t) &= dX^{u}(t)= b(t,X(t),Y(t), A(t), Z(t),K(t,\cdot), u(t),\omega)dt\nonumber\\
& + \sigma(t,X(t),Y(t),A(t), Z(t),K(t,\cdot),u(t),\omega)dB(t) \nonumber\\
&+ \int_{\RB} \gamma(t,X(t),Y(t),A(t), Z(t),K(t,\cdot),u(t),\zeta,\omega)\tN(dt,d\zeta)\; ; \; t \in [0,T] \nonumber \\
X(0)&=x \in \RB
\end{cases}
\end{equation}
\item{Predictive BSDE in $Y(t),Z(t),K(t)$:}
\begin{equation}\label{eq2.1}
\begin{cases}
dY(t)  &= -g(t, X(t),Y(t), A(t), Z(t),K(t,\cdot),u(t), \omega)dt + Z(t) dB(t)  \\
 & \quad+ \int_\RB K(t, \zeta) \tN (dt, d\zeta) \; ; \; t \in [0,T)  \\
%Y(t) & = L \; ; \; t \in [T, T+\delta] \; ; \; \delta > 0 \text{ fixed}.
Y(T) & = h(X(T),\omega). 
\end{cases}
\end{equation}
We set
\begin{equation}\label{eq2.2a}
Y(t):= L\; ; \; t \in (T, T + \delta] , 
\end{equation}
where $L$ is a given bounded $\mathcal{F}$-measurable random variable, representing a "cemetery" state of the process $Y$ after time $T$.
The process $A(t)$ represents our \emph{predictive mean-field} term. It is defined by 
\begin{equation}\label{eq2.2}
A(t) : = 
%\begin{cases}
E[Y(t + \delta) \mid \FC_t ] \; ; \; t \in [0, T]. 
%\end{cases}
\end{equation}
\end{itemize}
Here $\mathcal{R}$ is the set of functions from $\RB_0 := \RB \backslash \{0\}$ into $\RB$, 
$h(x,\omega)$ is a $C^1$ function (with respect to $x$) from $\RB \times \Omega$ into $\RB$ such that $h(x,\cdot)$ is $\mathcal{F}_T$-measurable for all $x$, 
and  
$$g : [0,T] \times \RB \times \RB \times \RB \times \RB \times \mathcal{R} \times \UB \times \Omega \rightarrow \RB$$
%, $\sigma : [0,T] \times \RB \times \RB \times \UB \times \Omega \rightarrow \RB$ and $\gamma : [0,T] \times \RB \times \RB \times \UB \times \RB_0 \times \Omega \rightarrow \RB$ 
\noindent
is a given function (driver) such that $g(t,x,y,a,z,k,u,\cdot)$
%, $\sigma(t,x,y,u,\cdot)$ and $\gamma(t,x,y,u,\zeta,\cdot)$ 
is an $\FB$-adapted process for all $x,y,a,z \in \RB, k\in\mathcal{R}$ and $u \in \UB$, which is the set of admissible control values. The process $u(t)$ is our control process, assumed to be in a given family $\AC = \AC_{\GB}$ of \emph{admissible} processes, assumed to be c\`adl\`ag and adapted to a given subfiltration $\GB=  \{ \GC_t\}_{t \geq 0}$ of the filtration $\FB$, i.e. $\mathcal{G}_t \subseteq \mathcal{F}_t$ for all $t$. The sigma-algebra $\mathcal{G}_t$ represents the information available to the controller at time $t$.

We assume that for all $u \in \mathcal{A}$ the coupled system \eqref{eq2.1a}-\eqref{eq2.2} has a unique solution $X(t)=X^{u}(t)\in L^{2}(m \times P),Y(t)=Y^{u}(t)\in L^{2}(m \times P), A(t)=A^{u}(t)\in L^{2}(m \times P),Z(t)=Z^{u}(t)\in L^{2}(m \times P),K(t,\zeta)=K^{u}(t,\zeta)\in L^{2}(m \times \nu \times P)$, with $X(t), Y(t), A(t)$ being c\`adl\`ag and $Z(t), K(t,\zeta)$ being predictable. Here and later $m$ denotes Lebesgue measure on $[0,T]$.\\

%We assume that $L$ is a given bounded $\FC_T$-measurable random variable representing the ``cemetery'' state of the process after time $T$.

To the best of our knowledge this system, \eqref{eq2.1a}-\eqref{eq2.2}, of \emph{predictive mean-field FBSDEs} has not been studied before.  However, the predictive BSDE \eqref{eq2.1}-\eqref{eq2.2} is related to the time-advanced BSDE which appears as an adjoint equation for stochastic control problems of a stochastic differential \emph{delay} equation. See  \cite{OSZ2} and the references therein. 

 The process $A(t)$ models the predicted future value of the state $Y$ at time $t + \delta$. Therefore \eqref{eq2.1}-\eqref{eq2.2} represent a system where the dynamics of the state is influenced by beliefs about the future. This is a natural model for situations where human behavior is involved, for example in pricing issues in  financial or energy markets . 
 
The performance functional associated to $u \in \AC$ is defined by
\begin{equation}\label{eq2.3}
J(u) = E \left[ \int_0^{T} f(t,X(t),Y(t), A(t), u(t),\omega)dt + \varphi(X(T),\omega)+\psi(Y(0)) \right]
\end{equation}
where $f : [0,T] \times \RB \times \RB \times \UB \times \Omega \rightarrow \RB$, $\varphi : \RB \times \Omega \rightarrow \RB$ and $\psi : \RB  \rightarrow \RB$ are given $C^{1}$ functions, with $f(t,x,y,a,u,\cdot)$ being $\FB$-adapted 
%and $\FC_T$-measurable
 for all $x,y,a \in \RB$, $u \in \UB$. We assume that $\varphi(x,\cdot)$ is $\mathcal{F}_{T}$-measurable for all $x$.

We study the following predictive mean-field stochastic control problem:

\paragraph{Problem} Find $u^* \in \AC$ such that
\begin{equation}\label{eq2.4}
\sup_{u \in \AC} J(u) = J(u^*).
\end{equation}
In Section \ref{sec3} we give a sufficient and a necessary maximum principle for the optimal control of forward-backward predictive mean-field systems of the type above.

An existence and uniqueness result for predictive mean-field BSDEs is given in Section \ref{sec4}. \\
Then in Section \ref{sec5} we apply the results to the following problems:
\begin{itemize}
\item
Portfolio optimization in a market where the stock price is modeled by a predictive mean-field BSDE,
\item
Optimization of consumption with respect to \emph{predictive recursive utility}.
\end{itemize}

\section{Solution methods for the stochastic control problem}\label{sec3}
\subsection{A sufficient maximum principle}\label{sec3.1}

For notational simplicity we suppress the dependence of $\omega$ in $f,g,h,\varphi$ and $\psi$ in the sequel.
We first give sufficient conditions for optimality of the control $u$  by modifying the stochastic maximum principle given in, for example, \cite{OS2}, to our new situation:\\
%Let $\RB$ denote the set of functions $r(\cdot) : \RB_0 \rightarrow \RB$, where $\RB_0 := \RB \backslash \{0\}$. 

We define the {\it Hamiltonian} $H : [0,T] \times \RB\times \RB\times \RB \times \RB \times\mathcal{R}\times \UB \times \RB\times \RB\times \RB\times \RB ) \rightarrow \RB$ associated to the problem \eqref{eq2.4} by
\begin{align}\label{eq3.1}
H(t,x,y,a,z,k,u,p,q,r,\lambda) & = f(t,x,y,a,u)+b(t,x,y,a,z,k,u)p+\sigma(t,x,y,a,z,k,u)q\nonumber\\
&+\int_{\RB}\gamma(t,x,y,a,z,k,u,\zeta)\tN(dt,d\zeta) + g(t,x,y,a,z,k,u)\lambda.
%+ b(t,x,y,u)p + \sigma(t,x,y,u)q \nonumber \\
% & + \int_\RB \gamma(t,x,y,u,\zeta) r(\zeta) \nu(d \zeta).
  \end{align}
% For simplicity of notation the dependence on $\omega$ is suppressed in the following.

We assume that $f,b,\sigma,\gamma$ and $g$, and hence $H$, are Fr\'echet differentiable $(C^1)$ in the variables $x,y,a,z,k,u$ and that the Fr\'echet derivative $\nabla_k H$ of $H$ with respect to $k \in \RC$ as a random measure is absolutely continuous with respect to $\nu$, with Radon-Nikodym derivative $\displaystyle \frac{d \nabla_k H}{d \nu}$. Thus, if $\langle \nabla_k H, h \rangle$ denotes the action of the linear operator $\nabla_k H$ on the function $h \in \RC$ we have
\begin{equation}
\langle \nabla_k H, h \rangle = \int_\RB h(\zeta) d \nabla_k H(\zeta) = \int_\RB h(\zeta) \frac{d \nabla_kH(\zeta)}{d \nu (\zeta)} d \nu (\zeta).
\end{equation}

%For $u \in \AC$ we let $(X^{u}(t),Y^u(t), Z^u(t), K^u(t,\cdot))$ be the associated solution of the coupled system \eqref{eq2.1a}-\eqref{eq2.2}, with $X(t)=X^{u}(t)$, $Y(t)=Y^{u}(t)$ being c\`adl\`ag, and $Z(t)=Z^{u}(t), K(t,\cdot)=K^{u}(t,\cdot)$ being predictable. We assume that for $u \in \AC$ these solutions exist and are unique and satisfy
%\begin{equation}\label{eq3.10}
%E \left[ \int_0^T \left\{ |X^{u}(t)|^2+|Y^u(t)|^2 + |Z^u(t)|^2 + \int_\RB | K^u(t,\zeta)|^2 \nu (d \zeta)\right\} dt \right] < \infty.
%\end{equation}

 The associated backward-forward system of equations in the {\it adjoint processes} $p(t),q(t),r(t),\lambda(t)$ is defined by
 \begin{itemize}
 \item{BSDE in $p(t),q(t),r(t)$:}
\begin{equation}\label{eq3.2a}
\begin{cases}
dp(t) & = \displaystyle - \frac{\partial H}{\partial x}(t)dt + q(t)dB(t) + \int_\RB r(t,\zeta) \tN (dt, d\zeta) \; ; \; 0 \leq t \leq T \\
p(T) & = \varphi'(X(T)) + \lambda(T) h'(X(T)).
\end{cases}
\end{equation}
\item{SDE in $\lambda(t)$:}
\begin{equation}\label{eq3.2}
\begin{cases}
d\lambda(t)  \displaystyle &=  \left\{ \frac{\partial H}{\partial y}(t) + \frac{\partial H}{\partial a}(t-\delta) \chi_{[\delta,T]}(t) \right\} dt + \frac{\partial H}{\partial z}(t) dB(t)  \\ & \qquad +  \displaystyle \int_\RB \frac{d \nabla_k H}{d \nu}(t,\zeta) \tN(dt, d\zeta) \; ; \; 0 \leq t \leq T \\
%+ q(t)dB(t)+ \int_\RB r(t,\zeta) \tN (dt, d\zeta) \; ; \; 0 \leq t \leq T \\
\lambda(0)& = \psi'(Y(0)),
\end{cases}
\end{equation}
\end{itemize}

where we have used the abbreviated notation
\begin{align*}
H(t) & = H(t, X(t),Y(t), A(t),Z(t),K(t,\cdot),u(t),p(t),q(t),r(t),\lambda(t)).
%\hH(t) & = H(t, \hX(t), \hY(t), \hu(t), \hp(t), \hq(t), \hr(t,\cdot))
\end{align*}
Note that, in contrast to the \emph{time advanced} BSDE \eqref{eq2.1}-\eqref{eq2.2}, \eqref{eq3.2} is a (forward) stochastic differential equation with \emph{delay}. 
%This type of BSDEs has been studied before (see \cite{DI2010} and the reference therein), but not in the context of optimal control of time-advanced mean-field equations.

\begin{theorem}(Sufficient maximum principle)\label{th3.1}
Let $\hu \in \AC$ with corresponding solution \\
$\hX(t),\hY(t), \hA(t), \hZ(t), \hK(t,\cdot), \hp(t),\hq(t),\hr(t),\hla(t)$ of \eqref{eq2.1a}-\eqref{eq2.2}, \eqref{eq3.2a}-\eqref{eq3.2}. Assume the following:
\begin{itemize}
\item
\begin{equation}\label{eq3.3a}
\hla(T) \geq 0  
\end{equation}
\item
 For all $t$, the functions
\begin{align}\label{eq3.3}
%\begin{cases}
&x \rightarrow h(x),x \rightarrow \varphi(x), x \rightarrow \psi(x) \text{ and }\nonumber \\
&(x,y,a,z,k,u) \rightarrow H(t,x,y,a,z,k,u,\hp(t),\hq(t),\hr(t),\hla(t)) \nonumber\\
\text{are concave}
%\end{cases}
\end{align}
\item
For all $t$ the following holds,
\begin{align}\label{eq3.15}
&\text{(The conditional maximum principle)} \nonumber\\
& \ess_{v \in \UB} E[H(t,\hX(t), \hY(t), \hA(t),\hZ(t), \hK(t,\cdot),v,\hla(t), \hp(t), \hq(t), \hr(t,\cdot)) \mid \GC_t] \nonumber\\
& \quad = E[H(t,\hX(t), \hY(t), \hA(t),\hZ(t), \hK(t,\cdot), \hu(t), \hla(t), \hp(t), \hq(t), \hr(t,\cdot)) \mid \GC_t] \; ; \; t \in [0,T]
\end{align}
\item
\begin{equation}\label{eq3.16}
\left\| \frac{d \nabla_k \hH(t,.)}{d \nu} \right\| < \infty \text{ for all } t \in [0,T].
\end{equation}
\end{itemize}
Then $\hu$ is an optimal control for the problem \eqref{eq2.4}.
\end{theorem}
\dproof By replacing the terminal time $T$ by an increasing sequence of stopping times $\tau_n$ converging to $T$ as $n$ goes to infinity, and arguing as in \cite{OS2} we see that we may assume that all the local martingales appearing in the calculations below are martingales. 

Much of the proof is similar to the proof of Theorem 3.1 in \cite{OS2}, but due to the predictive mean-field feature of the BSDE \eqref{eq2.1}-\eqref{eq2.2}, there are also essential differences. Therefore, for the convenience of the reader, we sketch the whole proof:

Choose $u \in \AC$ and consider
\begin{equation}\label{eq3.5}
J(u) - J(\hu) = I_1 + I_2+I_3,
\end{equation}
with
\begin{equation}\label{eq3.6}
I_1 := E \left[ \int_0^{T} \{f(t) - \hf(t)\} \right]dt,\quad I_2 := E[\varphi(X(T)) - \varphi(\hX(T))],\quad I_3 := \psi(Y(0)) - \psi(\hY(0)),
\end{equation}
where $\hf(t) = f(t, \hY(t), \hA(t), \hu(t))$ etc., and $\hY(t) = Y^{\hu}(t)$ is the solution of \eqref{eq2.1}-\eqref{eq2.2} when $u = \hu$, and $\hA(t) = E[ \hY(t) \mid \FC_t]$.

By the definition of $H$ we have
\begin{equation}\label{eq3.7}
I_1 = E \left[ \int_0^{T} \{H(t) - \hH(t) -\hp(t)\tb(t)
- \hq(t) \tsi(t) - \int_\RB \hr(t,\zeta) \tga(t,\zeta) \nu(d\zeta)- \hla(t) \tg(t) 
\right],
\end{equation}
where we from now on use the abbreviated notation
\begin{align*}
H(t) & = H(t, X(t),Y(t), A(t), Z(t),K(t,\cdot),u(t), \hla(t)) \\
\hH(t) & = H(t, \hX(t),\hY(t), \hA(t),\hZ(t),\hK(t,\cdot), \hu(t), \hla(t))
\end{align*}
and we put
$$\tb(t) := b(t) - \hb(t),$$
and similarly with $\tX(t):=X(t)-\hX(t), \tY(t):=Y(t)-\hY(t),\tA(t):=A(t)-\hA(t), $ etc.

By concavity of $\varphi$, \eqref{eq3.2} and the It\^o formula,
\begin{align}\label{eq3.13}
I_2 & \leq E[ \varphi'(\hX(T))\tilde{X}(T)] \nonumber \\
& = E [ \hp(T) \tX(T)] - E [ \hla(T) h'(\hX(T))\tilde{X}(T)] \nonumber \\
& = \left( E \left[ \int_0^{T} \hp(t^-) d\tX(t) + \int_0^{T} \tX(t^-) d \hp (t) \right. \right. \nonumber 
+ \int_0^{T} \hq(t) \tilde{\sigma}(t)dt \nonumber \\
& \left. \left. + \int_0^{T} \int_\RB \hr(t,\zeta)\tilde{\gamma}(t,\zeta) \nu (d \zeta)dt \right]\right) - E [ \hla(T) h'(\hX(T))\tX(T)] \nonumber \\
& = E \left[ \int_0^{T} \hp(t) \tilde{b}(t)dt + \int_0^{T} \tX(t) \left( - \frac{\partial \hH}{\partial x}(t)\right)dt \right.\nonumber \\
& \left.\quad + \int_0^{T} \hq(t)\tilde{\sigma}(t)dt + \int_0^{T}\int_\RB \hr(t,\zeta) \tilde{\gamma}(t,\zeta) \nu(d\zeta) dt \right] \nonumber \\
& \quad - E [ \hla(T) h'(\hX(T))\tilde{X}(T)].
\end{align}

By concavity of $\psi$ and $h$, \eqref{eq3.3a} and the It\^o formula we have
\begin{align}\label{eq3.14a}
I_3 & \leq E \left[  \psi'(Y(0)) \tY(0)\right] = E[\hla(0) \tY(0)] \nonumber \\
& = E[\hla(T) \tY(T)] - E \left[ \int_0^T \hla(t) d\tY(t) + \int_0^T \tY(t) d \hla(t)+\int_0^T d[\tY,\hla](t) \right] \nonumber\\
& = E[\hla(T) (h(X(T)) - h(\hX(T)))] - E \left[ \int_0^T \hla(t) d\tY(t) + \int_0^T \tY(t) d \hla(t)+\int_0^T d[\tY,\hla](t) \right] \nonumber\\
& \leq E[\hla(T) h'(\hX(T))\tX(T)] + E \left[ \int_0^T \hla(t) \tg(t) dt + \int_0^T \tY(t) \left[ - \frac{\partial \hH}{\partial y}(t) - \frac{\partial \hH}{\partial a} (t-\delta) \chi_{[\delta,T]}(t)\right] dt \right. \nonumber \\
 & \quad + \int_0^T \frac{\partial \hH}{\partial z}(t) \tZ(t)dt 
  \quad \left. + \int_0^T \int_\RB \frac{d \nabla_k \hH}{d \nu}(t,\zeta) \tK(t,\zeta) \nu(d\zeta)dt 
 \right]
%& \left. \quad + \int_0^T \hq(t) \tsi(t) dt + \int_0^T \int_\RB \hr(t,\zeta) \tga(t,\zeta) \nu(d\zeta)dt \right]
\end{align}

Adding \eqref{eq3.7},\eqref{eq3.13} and \eqref{eq3.14a} we get, by \eqref{eq3.2},
\begin{align}\label{eq3.11}
J(u) & - J(\hu) = I_1 +I_2+I_3  \nonumber \\
& \leq E \left[ \int_0^T\left\{ H(t) - \hH(t) - \frac{\partial \hH}{\partial x} \tX(t)- \frac{\partial \hH}{\partial y} \tY(t) \nonumber \right. \right. \\
& \quad - \frac{\partial H}{\partial a}(t-\delta) \chi_{[\delta,T]}(t) \tY(t) - \frac{\partial H}{\partial z}(t) \tZ(t)
 \left.  - \langle \nabla_k \hH(t,\cdot), \tK(t,\cdot)\rangle  \biggr{\}} dt \right].
\end{align}

Note that, since $Y(s) = \hY(s) = L$ for $s \in (T,T+\delta]$ by \eqref{eq2.1a}, we get
\begin{align}\label{eq3.12}
E &\left[ \int_0^T \frac{\partial \hH}{\partial a} (t- \delta) \tY(t) \chi_{[\delta,T]}(t) dt \right] = E \left[ \int_0^{T-\delta} \frac{\partial \hH}{\partial a}(s) \tY(s + \delta)ds \right] \nonumber \\
& = E \left[ \int_0^{T - \delta} E \left[ \frac{\partial \hH}{\partial a}(s) \tY(s + \delta) \mid \FC_s\right]dt \right] \nonumber \\
& = E \left[ \int_0^{T-\delta} \frac{\partial \hH}{\partial a}(s) E \left[ \tY(s + \delta) \mid \FC_s\right] ds \right] = E \left[ \int_0^T \frac{\partial \hH}{\partial a}(s) \tA(s)ds \right].
\end{align}

Substituted into \eqref{eq3.11} this gives, by concavity of $H$,
\begin{align}\label{eq3.22a}
J(u) & - J(\hu) = I_1 +I_2 +I_3\nonumber \\
& \leq E \left[ \int_0^T\left\{ H(t) - \hH(t) - \frac{\partial \hH}{\partial x} (X(t) - \hX(t))- \frac{\partial \hH}{\partial y} (Y(t) - \hY(t)) \nonumber \right. \right. \\
& \quad - \frac{\partial H}{\partial a}(t) (A(t) - \hA(t)) - \frac{\partial H}{\partial z}(t) (Z(t) - \hZ(t)) \nonumber \\
& \quad \left.  - \langle \nabla_k \hH(t,\cdot), (K(t,\cdot) - \hK(t,\cdot)\rangle  \biggr{\}} dt \right]\nonumber\\
& \leq E \left[ \int_0^T \frac{\partial \hH}{\partial u} (t) (u(t) - \hu(t)) dt\right]\nonumber\\
& = E \left[ \int_0^T E[\frac{\partial \hH}{\partial u} (t)|\mathcal{G}_t] (u(t) - \hu(t)) dt  \right] \leq 0,
\end{align}
since $u = \hu(t)$ maximizes $E[\hH(t)|\mathcal{G}_t]$.
\fproof

\subsection{A necessary maximum principle}\label{sec3.2}

We proceed to prove a partial converse of Theorem \ref{th3.1}, in the sense that we give {\it necessary} conditions for a control $\hu$ to be optimal. In this case we can only conclude that $\hu(t)$ is a critical point for the Hamiltonian, not necessarily a maximum point. On the other hand, we do not need any concavity assumptions, but instead we need some properties of the set $\AC$ of admissible controls, as described below.

%In the following $m$ denotes Lebesgue measure on $[0,T]$.

\begin{theorem}(Necessary maximum principle)\label{th3.2}

Suppose $\hu \in \AC$ with associated solutions $\hX, \hY, \hZ, \hK, \hat{p},\hat{q},\hat{r},\hla$ of \eqref{eq2.1a}-\eqref{eq2.2} and \eqref{eq3.2a}-\eqref{eq3.2}. Suppose that for all processes $\beta(t)$ of the form
\begin{equation}\label{eq3.14}
\beta(t) := \chi_{[t_0, T]}(t) \alpha,
\end{equation}
where $t_0 \in [0,T)$ and $ \alpha = \alpha(\omega)$ is a bounded $\GC_{t_0}$-measurable random variable, there exists $\delta > 0$ such that the process
$$ \hu(t) + r\beta(t) \in \AC \text{ for all } r \in [- \delta, \delta].$$

We assume that the {\it derivative processes} 
%$y(t)$ of $Y^u(t)$ at $u = \hu$, 
defined by
\begin{equation}%\label{eq3.18}
x(t) = x^\beta(t) = \frac{d}{dr} X^{\hu + r \beta}(t) \mid_{r=0},
\end{equation}
\begin{equation}%\label{eq3.19}
y(t) = y^\beta(t) = \frac{d}{dr} Y^{\hu + r \beta}(t) \mid_{r=0},
\end{equation}
\begin{equation}%\label{eq3.20}
a(t) = a^\beta(t) = \frac{d}{dr} A^{\hu + r \beta}(t) \mid_{r=0},
\end{equation}
\begin{equation}%\label{eq3.21}
z(t) = z^\beta(t) = \frac{d}{dr} Z^{\hu + r \beta}(t) \mid_{r=0},
\end{equation}
\begin{equation}%\label{eq3.22}
k(t) = k^\beta(t) = \frac{d}{dr} K^{\hu + r \beta}(t) \mid_{r=0},
\end{equation}
exist and belong to  $L^2(m \times P)$, $L^2(m \times P)$, $L^2(m \times P)$, and $L^2(m \times P \times \nu)$, respectively.

Moreover, we assume that $x(t)$ satisfies the equation
\begin{equation}\label{eq3.23}
\begin{cases}
dx(t) =  \displaystyle \left\{ \frac{\partial b}{\partial x}(t) x(t) + \frac{\partial b}{\partial y}(t)y(t)+ \frac{\partial b}{\partial a}(t)a(t) + \frac{\partial b}{\partial z}(t)z(t) + \langle\nabla_k b, k(t,\cdot)\rangle  + \frac{\partial b}{\partial u}(t) \beta(t) \right\} dt \\
 \quad \displaystyle + \left\{ \frac{\partial \sigma}{\partial x}(t) x(t) + \frac{\partial \sigma}{\partial y}(t) y(t) + \frac{\partial \sigma}{\partial a}(t) a(t)+ \frac{\partial \sigma}{\partial z}(t) z(t) + \langle \nabla_k \sigma, k(t,\cdot)\rangle + \frac{\partial \sigma}{\partial u}(t) \beta(t)\right\} dB(t) \\
 \quad \displaystyle + \int_\RB
 \left\{ \frac{\partial \gamma}{\partial x}(t,\zeta) x(t) + \frac{\partial \gamma}{\partial y}(t,\zeta) y(t)+ \frac{\partial \gamma}{\partial a}(t,\zeta) a(t) + \frac{\partial \gamma}{\partial z}(t,\zeta) z(t) + \langle \nabla_k \gamma(t,\zeta), k(t,\cdot)\rangle \right.\\
 \qquad \left. \displaystyle + \frac{\partial \gamma}{\partial u}(t,\zeta) \beta(t)\right\}  \tN(dt,d\zeta) \; ; \; t \in [0,T] \\
 x(0) = 0
\end{cases}
\end{equation}
and that $y(t)$ satisfies the equation
\begin{align}\label{eq3.25}
\begin{cases}
dy(t)  &= -\left\{\frac{\partial g}{\partial x}(t) x(t)+\frac{\partial g}{\partial y}(t) y(t) + \frac{\partial g}{\partial a}(t)a(t) + \frac{\partial g}{\partial z}(t)z(t) + \langle \nabla_k g(t), k(t,\cdot)\rangle + \frac{\partial g}{\partial u}(t) \beta(t)   \right\} dt  \\
& + z(t) dB(t) + \int_\RB k(t,\zeta) \tN(dt, d\zeta)\; ; \; 0 \leq t < T \\
y(T)&=h'(X(T))x(T)\\
y(t)&= 0 \; ; \; T < t \leq T+\delta,
\end{cases}
\end{align}where we have used the abbreviated notation
\begin{equation}
\frac{\partial g}{\partial x}(t) = \frac{\partial}{\partial x} g(t,x,y,a,z,k,u)_{x=X(t),y =Y(t),a=A(t),z=Z(t),k=K(t),u=u(t)}  \text{ etc.} \nonumber
\end{equation}

Then the following, (i) and (ii), are equivalent:

\begin{myenumerate}
\item $\displaystyle \frac{d}{dr} J(\hu + r \beta)_{r=0} = 0$ for all $\beta$ of the form \eqref{eq3.14}
\item $\displaystyle \frac{d}{du} E[H(t, \hY(t), \hA(t),\hZ(t),\hK(t),u, \hla(t))_{u = \hu(t)}|\mathcal{G}_t] = 0,$ \\
where $(\hY, \hA, \hZ,\hK,\hla)$ is the solution of \eqref{eq2.1},\eqref{eq2.2} and \eqref{eq3.2} corresponding to $u=\hu$.
\end{myenumerate}
\end{theorem}

\dproof
As in Theorem 3.1, by replacing the terminal time $T$ by an increasing sequence of stopping times $\tau_n$ converging to $T$ as n goes to infinity, we obtain as in \cite{OS2} that we may assume that all the local martingales appearing in the calculations below are martingales. The proof has many similarities with the proof of Theorem 3.2 in \cite{OS2}, but since there are some essential differences due to the predictive mean-field term, we sketch the whole proof. For simplicity of notation we drop the hats in the sequel, i.e. we write $u$ instead of $\hu$ etc.\\
\noindent (i) $\Rightarrow$ (ii): 
We can write $\displaystyle \frac{d}{dr} J(u + r \beta) \mid_{r=0} = I_1 + I_2 +I_3$, where
\begin{align*}
I_1 & = \frac{d}{dr} E \left[ \int_0^T f(t,Y^{u+r \beta}(t), A^{u+r \beta}(t), Z^{u+r \beta}(t), K^{u+r \beta}(t), u(t) + r \beta(t))dt\right]_{r=0} \\
I_2 & = \frac{d}{dr} [ \varphi(X^{u+r \beta}(T))]_{r=0}\\
I_3 & = \frac{d}{dr} [ \psi(Y^{u+r \beta}(0))]_{r=0}.
\end{align*}
By our assumptions on $f$ and $\psi$ we have
\begin{align}\label{eq3.26}
I_1&= \left[ \int_0^T \left\{ \frac{\partial f}{\partial x}(t) x(t) +\frac{\partial f}{\partial y}(t) y(t) + \frac{\partial f}{\partial a}(t) a(t) + \frac{\partial f}{\partial z}(t) z(t) + \langle\nabla_k f(t,\cdot), k(t,\cdot)\rangle + \frac{\partial f}{\partial u}(t) \beta(t)\right\} dt\right]\\
 \text{ and }  \nonumber\\
 \label{eq3.27a}
 I_2&=E[\varphi'(X(T)x(T)]=E[p(T)x(T)]\\
\label{eq3.28}
I_3 &  = \psi'(Y(0))y(0) = \lambda(0) y(0).
\end{align}
By the It\^o formula and \eqref{eq3.23}
\begin{align}\label{eq3.27}
I_2 & = E[p(T) x(T)]  \nonumber = E \left[ \int_0^{T} p(t) dx(t) + \int_0^{T} x(t) dp(t) + \int_0^{T} d[p,x](t) \right] \nonumber \\
& = E \left[ \int_0^{T} p(t) \left\{ \frac{\partial b}{\partial x}(t) x(t) + \frac{\partial b}{\partial y}(t) y(t)+ \frac{\partial b}{\partial a}(t) a(t) + \frac{\partial b}{\partial z}(t) z(t) + \langle \nabla_k  b(t),k(t,\cdot)\rangle \right.\right.\nonumber \\
& \left. \quad + \frac{\partial b}{\partial u}(t) \beta(t) \right\} dt + \int_0^{\tau_n}x(t) \left( - \frac{\partial H}{\partial x}(t)\right) dt  + \int_0^{\tau_n} q(t) \left\{ \frac{\partial \sigma}{\partial x}(t) x(t)  \right. \nonumber \\
& \quad \left. + \frac{\partial \sigma}{\partial y}(t) y(t)+ \frac{\partial \sigma}{\partial a}(t) a(t)  + \frac{\partial \sigma}{\partial z}(t) z(t)  + \langle \nabla_k \sigma(t), k(t,\cdot)\rangle + \frac{\partial \sigma}{\partial u}(t) \beta(t) \right\}dt \nonumber \\
& \quad + \int_0^{T}\int_\RB r(t,\zeta) \left\{ \frac{\partial \gamma}{\partial x}(t,\zeta)x(t) + \frac{\partial \gamma}{\partial y}(t,\zeta)y(t)+ \frac{\partial \gamma}{\partial a}(t,\zeta)a(t) + \frac{\partial \gamma}{\partial z}(t,\zeta)z(t) + < \nabla_k \gamma(t,\zeta), k(t,\cdot)> \right.\nonumber \\
& \quad \left. \left. + \frac{\partial \gamma}{\partial u}(t,\zeta)\beta(t)\right\} \nu (d\zeta) dt \right] \nonumber \\
& = E \left[ \int_0^{T} x(t) \left\{ \frac{\partial b}{\partial x}(t) p(t) + \frac{\partial \sigma}{\partial x}(t) q(t) + \int_\RB \frac{\partial \gamma}{\partial x}(t,\zeta) r(t,\zeta) \nu(d\zeta) - \frac{\partial H}{\partial x}(t) \right\} dt\right. \nonumber \\
& \quad + \int_0^{T} y(t) \left\{ \frac{\partial b}{\partial y}(t) p(t) + \frac{\partial \sigma}{\partial y}(t) q(t) + \int_\RB \frac{\partial \gamma}{\partial y}(t,\zeta)  r(t,\zeta) \nu(d\zeta)  \right\} dt\nonumber \\
& \quad + \int_0^{T} a(t) \left\{ \frac{\partial b}{\partial a}(t) p(t) + \frac{\partial \sigma}{\partial a}(t) q(t) + \int_\RB \frac{\partial \gamma}{\partial a}(t,\zeta)  r(t,\zeta) \nu(d\zeta)  \right\} dt\nonumber \\
& \quad + \int_0^{T} z(t) \left\{ \frac{\partial b}{\partial z}(t) p(t) + \frac{\partial \sigma}{\partial z}(t) q(t) + \int_\RB \frac{\partial \gamma}{\partial z}(t,\zeta)  r(t,\zeta) \nu(d\zeta)  \right\} dt\nonumber \\
& \quad \left. + \int_0^{T}\int_\RB   \langle k(t,\cdot),  \nabla_kb(t)p(t) + \nabla_k \sigma(t) q(t) + \int_\RB \nabla_k \gamma(t,\zeta) r(t,\zeta) \nu (d \zeta) \rangle \nu(d\zeta)dt\right]\nonumber \\
& = E \left[ \int_0^{T} x(t) \left\{ - \frac{\partial f}{\partial x}(t) - \lambda (t) \frac{\partial g}{\partial x}(t)  \right\} dt + \int_0^{T} y(t) \left\{ \frac{\partial H}{\partial y}(t)- \frac{\partial f}{\partial y}(t) - \lambda(t)  \frac{\partial g}{\partial y}(t)    \right\} dt \right. \nonumber \\
& \quad + \int_0^{T} a(t) \left\{ \frac{\partial H}{\partial a}(t)- \frac{\partial f}{\partial a}(t) - \lambda(t)  \frac{\partial g}{\partial a}(t)    \right\} dt+ \int_0^{T} z(t) \left\{ \frac{\partial H}{\partial z}(t)- \frac{\partial f}{\partial z}(t) - \lambda(t)  \frac{\partial g}{\partial z}(t)    \right\} dt\nonumber \\
& \quad + \int_0^{T} \int_\RB k(t,\zeta) \{ \nabla_kH(t) - \nabla_k f(t) - \lambda(t) \nabla_k g(t)\} \nu (d\zeta) dt \nonumber \\
& \quad \left. + \int_0^{T} \beta(t) \left\{ \frac{\partial H}{\partial u}(t)- \frac{\partial f}{\partial u}(t) - \lambda(t)  \frac{\partial g}{\partial u}(t)    \right\} dt \right]\nonumber \\
& = - I_1 - E \left[ \int_0^T \lambda(t) \left\{ \frac{\partial g}{\partial x}(t) x(t) + \frac{\partial g}{\partial y} (t) y(t) + \frac{\partial g}{\partial z} (t) z(t) \right. \right.\nonumber \\
& \quad \left. + \langle \nabla_kg(t), k(t,\cdot)\rangle + \frac{\partial g}{\partial u}(t) \beta(t) \right\} dt \nonumber \\
& \quad + E \left[ \int_0^T \left\{ \frac{\partial H}{\partial y} (t) y(t) + \frac{\partial H}{\partial z} (t) z(t)  + \langle \nabla_k H(t), k(t,\cdot)\rangle + \frac{\partial H}{\partial u} (t) \beta(t)\right\}dt\right]
\end{align}
By the It\^o formula and \eqref{eq3.25},
\begin{align}\label{eq3.22}
I_3 & = \lambda(0) y(0) = 
E \left[ \lambda(T) y(T)
  -  \left( \int_0^{T} \lambda(t) dy(t) + \int_0^{T} y(t) d \lambda(t) + \int_0^{T} d[\lambda, y](t)\right)\right] \nonumber \\
& = E[ \lambda(T) y(T)] \nonumber \\
&\quad - \left( E \left[ \int_0^{T} \lambda(t) \left\{ - \frac{\partial g}{\partial y}(t) y(t) - \frac{\partial g}{\partial a}(t) a(t) - \frac{\partial g}{\partial z}(t) z(t) \right. \right. \right. \nonumber \\
& \left. \quad - \langle \nabla_k g(t), k(t,\cdot)\rangle - \frac{\partial g}{\partial u}(t) \beta(t) \right\}dt \nonumber \\
& + \int_0^{T} y(t) \frac{\partial H}{\partial y}(t) dt + y(t) \frac{\partial H}{\partial a}(t-\delta) \chi_{[\delta,T]}(t) dt 
+\int_0^{T} z(t) \frac{\partial H}{\partial z}(t) dt \nonumber \\
& \left. \left.\quad 
+
\int_0^{T}   \int_\RB k(t,\zeta) \nabla_k H(t,\zeta) \nu (d \zeta) dt \right] \right).
\end{align}
Adding \eqref{eq3.26}, \eqref{eq3.27} and \eqref{eq3.22} and using that
\begin{align}\label{eq3.17}
&E[ \int_0^T y(t)\frac{\partial H}{\partial a}(t-\delta) \chi_{[\delta,T]} dt ]
=E[ \int_0^{T-\delta} y(s+\delta)\frac{\partial H}{\partial a}(s) ds]\nonumber\\
&=E[ \int_0^T \frac{\partial H}{\partial a}(s)E[y(s+\delta)|\mathcal{F}_s] ds] 
=E[ \int_0^T y(t)\frac{\partial H}{\partial a}(s) a(s) ds],  
\end{align}
we get
$$ \frac{d}{dr} J(u + r \beta) \mid_{r=0} = I_1 + I_2 = E \left[  \int_0^T \frac{\partial H}{\partial u}(t) \beta(t)dt \right].$$
We conclude that
$$\frac{d}{dr} J(\hu + r \beta) \mid_{r=0}  = 0  $$
if and only if
$$E \left[ \int_0^T \frac{\partial \hH}{\partial u}(t) \beta(t) dt \right] = 0 \;  \; \text{ for all bounded } \beta \in \AC_{\mathbb{G}} \text { of the form } \eqref{eq3.14}.$$
Since this holds for all such $\beta$, we obtain that if (i) holds, then
\begin{equation}\label{eq3.24}
\int_{t_0}^T E \left[ \frac{\partial \hH}{\partial u}(t) \mid \GC_{t_0} \right]dt = 0 \text{ for all } t_0 \in [0,T).
\end{equation}
Differentiating with respect to $t_0$ and using continuity of $\displaystyle \frac{\partial \hH}{\partial u}(t)$, we conclude that (ii) holds.

\medskip

\noindent (ii) $\Rightarrow$ (i): This is proved by reversing the above argument. We omit the details.
\fproof

\section{Existence and uniqueness of predictive mean-field equations} \label{sec4}
In this section we study the existence and uniqueness of predictive mean-field BSDEs in the unknowns $Y(t),Z(t),K(t,\zeta)$ of the form
\begin{align}\label{eq4.1}
\begin{cases}
dY(t) & = -g(t,Y(t), A(t), Z(t),K(t,\cdot), \omega)dt + Z(t) dB(t)  \\
& + \int_\RB K(t, \zeta) \tN (dt, d\zeta) \; ; \; t \in [0,T)  \\
Y(t) & = L \; ; \; t \in [T, T+\delta] \; ; \; \delta > 0 \text{ fixed},
\end{cases}
\end{align}
\noindent where $L \in L^2(P)$ is a given $\mathcal{F}_T$-measurable random variable, and the process $A(t)$ as before is defined by 
\begin{equation}\label{eq4.2}
A(t) = 
%\begin{cases}
E[Y(t + \delta) \mid \FC_t ] \; ; \; t \in [0, T].
%\end{cases}
\end{equation}
To this end, we can use the same argument which was used to handle a similar, but different, time-advanced BSDE in \cite{OSZ2}. For completeness we give the details:

\begin{theorem}
Suppose the following holds
\begin{align}
%\begin{itemize}
%\item
&E[\int_0^T g^2(t,0,0,0,0)dt] < \infty \\
%\item
&\text{There exists a constant } C \text{ such that }\nonumber\\
&|g(t,y_1,a_1,z_1,k_1) -g(t,y_2,a_2,z_2,k_2)|\leq C(|y_1-y_2| + |z_1-z_2| +
(\int_{\mathbb{R}}|k_1(\zeta) - k_2(\zeta)|^2\nu(d\zeta))^{\frac{1}{2}})
%\end{itemize}
\end{align}
for all $t \in [0,T],$ a.s.
Then there exists a unique solution triple $(Y(t),Z(t),K(t,\zeta))$ of \eqref{eq4.1} such that the following holds:
\begin{align}
\begin{cases}
&Y \text{ is c\` adl\` ag and } E[\sup_{t \in [0,T]} Y^2(t)] < \infty, \nonumber\\
&Z,K \text{ are predictable and } E[\int_{0}^T\{ Z^2(t) + \int_{\mathbb{R}}K^2(t,\zeta)\nu(d\zeta)\}dt] < \infty.
\end{cases}
\end{align}

\end{theorem}

\dproof
We argue backwards, starting with the interval $[T-\delta,T]$:\\
\emph{Step 1.}  In this interval we have $A(t)=E[L|\mathcal{F}_t]$ and hence we know from the theory of classical BSDEs (see e.g. \cite{Q},\cite{QS} and the references therein), that there exists a unique solution triple $(Y(t),Z(t),K(t,\zeta))$ such that the following holds:
\begin{align}
\begin{cases}
&Y \text{ is c\` adl\` ag and } E[\sup_{t \in [T-\delta,T]} Y^2(t)] < \infty, \nonumber\\
&Z,K \text{ are predictable and } E[\int_{T-\delta}^T\{ Z^2(t) + \int_{\mathbb{R}}K^2(t,\zeta)\nu(d\zeta)\}dt] < \infty.
\end{cases}
\end{align}
\emph{Step 2.} Next, we continue with the interval $[T-2\delta,T-\delta]$. For $t$ in this interval, the value of $Y(t+\delta)$ is known from the previous step and hence $A(t)=E[Y(t+\delta)|\mathcal{F}_t]$ is known. Moreover, by Step1 the terminal value for this interval, $Y(T-\delta)$, is known and in $L^2(P)$. Hence we can again refer to the theory of classical BSDEs and get a unique solution in this interval.\\
\emph{Step n}. We continue this iteration until we have reached the interval $[0, T-n\delta]$, where $n$ is a natural number such that 
$$T-(n+1)\delta \leq 0 < T-n\delta.$$
Combining the solutions from each of the subintervals, we get a solution for the whole interval.
\fproof

\section{Applications}\label{sec5}
In this section we illustrate the results of the previous sections by looking at two examples.
\subsection{Optimal portfolio in an insider influenced market}
In the seminal papers by Kyle \cite{K} and Back \cite{B} it is proved that in a financial market consisting of 
\begin{itemize}
\item
\emph{noise traders} (where noise is modeled by Brownian motion), 
\item
\emph{an insider} who knows the value $L$ of the price of the risky asset at the terminal time $t=T$ and 
\item
\emph{a market maker} who at any time $t$ clears the market and sets the market price, 
\end{itemize}
the corresponding equilibrium price process (resulting from the insider's portfolio which maximizes her expected profit), will be a \emph{Brownian bridge terminating at the value $L$ at time $t =T$}. 
In view of this we see that a predictive mean-field equation can be a natural model of the \emph{risky asset price in an insider influenced market}.\\

Accordingly, suppose we have a market with the following two investment possibilities:
\begin{itemize}
\item
A risk free asset, with unit price $S_0(t):=1$ for all $t$
\item
A risky asset with unit price $S(t):=Y(t)$ at time $t$, given by the predictive mean-field equation
\begin{align}\label{eq5.1}
\begin{cases}
dY(t) & = -A(t)\mu(t)dt + Z(t) dB(t); \; t \in [0,T) \\
Y(t) & = L(\omega); \quad t \in [T,T+\delta],
\end{cases}
\end{align}
\noindent where $\mu(t)=\mu(t,\omega)$ is a given bounded adapted process and $L$ is a given bounded $\mathcal{F}_T$-measurable random variable, being the terminal state of the process $Y$ at time $T$.
\end{itemize}
Let $u(t)$ be a portfolio, representing the number of risky assets held at time $t$. We assume that $\mathbb{G} = \mathbb{F}$.
If we assume that the portfolio is self-financing, the corresponding wealth process $X(t)=X^{u}(t)$ is given by
\begin{equation}\label{eq5.2}
\begin{cases}
dX(t)= u(t) dY(t) = u(t) A(t)\mu(t) dt + u(t) Z(t) dB(t); \; t \in [0,T)\\
X(0) = x > 0.
\end{cases}
\end{equation}

Let $U: [0,\infty) \mapsto [-\infty,\infty)$ be a given utility function, assumed to be increasing, concave and $C^{1}$ on $(0,\infty)$. We study the following portfolio optimization problem:

\paragraph{Problem 5.1} Find $u^* \in \AC$ such that
\begin{equation}\label{eq5.3}
\sup_{u \in \AC} E[U(X^{u}(T))] = E[U(X^{u^*}(T))].
\end{equation}

This is a problem of the type discussed in the previous sections, with $f=\psi=N=0, \varphi=U$ and $ h(x,\omega)=L(\omega)$, and we can apply the maximum principles from Section 3 to study it. 

By \eqref{eq3.1} the Hamiltonian gets the form
\begin{align}\label{eq5.4}
H(t,x,y,a,z,k,u,p,q,r,\lambda) & = ua\mu(t)p+uzq+a\mu(t)\lambda.
  \end{align}
  
 The associated backward-forward system of equations in the adjoint processes $p(t),q(t),\lambda(t)$ becomes
 \begin{itemize}
 \item{BSDE in $p(t),q(t)$:}
\begin{equation}\label{eq5.5}
\begin{cases}
dp(t) & =  q(t)dB(t) \; ; \; 0 \leq t \leq T \\
p(T) & = U'(X(T)),
\end{cases}
\end{equation}
\item{SDE in $\lambda(t)$:}
\begin{equation}\label{eq5.6}
\begin{cases}
d\lambda(t) =  \mu(t-\delta)[u(t-\delta)p(t-\delta)+\lambda(t-\delta)] \chi_{[\delta,T]}(t) dt + u(t)q(t) dB(t)\; ; \; 0 \leq t \leq T \\
\lambda(0) = 0.
\end{cases}
\end{equation}
\end{itemize}  
  
 The Hamiltonian can only have a maximum with respect $u$ if 
 \begin{equation}
 A(t)\mu(t)p(t)+Z(t)q(t)=0.
 \end{equation}
Substituting this into \eqref{eq5.5} we get
\begin{equation}\label{eq5.8}
\begin{cases}
dp(t) & = -\theta(t) p(t)dB(t) \; ; \; 0 \leq t \leq T \\
p(T) & = U'(X(T)),
\end{cases}
\end{equation}
where
\begin{equation}\label{eq5.9}
\theta(t):= \frac{A(t)\mu(t)}{Z(t)}.
\end{equation}
From this we get
\begin{equation}\label{eq5.10}
p(t)=c \exp(-\int_0^t \theta(s)dB(s) -\frac{1}{2}\int_0^t (\theta(s))^2 ds)\; ; \; 0 \leq t \leq T
\end{equation}
where the constant 
\begin{equation}\label{eq5.11}
c =p(0)= E[U'(X(T)] 
\end{equation}
remains to be determined.

In particular, putting $t=T$ in \eqref{eq5.10} we get
\begin{equation}\label{eq5.12}
U'(X(T))= p(T)=c \exp(-\int_0^T \theta(s)dB(s) -\frac{1}{2}\int_0^T (\theta(s))^2 ds)
\end{equation}
or
\begin{equation}\label{eq5.13}
X(T)= (U')^{-1}(c \exp(-\int_0^T\theta(s)dB(s) -\frac{1}{2}\int_0^T (\theta(s))^2 ds)).
\end{equation}
Define
\begin{equation}\label{eq5.14}
\Gamma(T)= \exp(\int_0^T \theta(s)dB(s) -\frac{1}{2}\int_0^T (\theta(s))^2 ds).
\end{equation}
Then by the Girsanov theorem the measure $Q$ defined on $\mathcal{F}_T$ by
\begin{equation}\label{eq5.15}
dQ(\omega)=\Gamma(T) dP(\omega)
\end{equation}
is an equivalent martingale measure for the market \eqref{eq5.1}.Therefore, by \eqref{eq5.13},
\begin{equation}\label{eq5.16}
x= E_Q[X(T)] =E[(U')^{-1}(c \exp(-\int_0^T\theta(s)dB(s) -\frac{1}{2}\int_0^T (\theta(s))^2 ds))\Gamma(T)].
\end{equation}
This equation determines implicitly the value of the constant $c$ and hence by \eqref{eq5.13} the optimal terminal wealth $X(T)=X^{u^*}(T)$. To find the corresponding optimal portfolio $u^*$ we proceed as follows:\\

Define 
\begin{equation}%\label{eq5.17}
Z_0(t) := u^*(t)Z(t). 
\end{equation}
Then $(X^{u^*}(t),Z_0(t))$ is found by solving the linear BSDE
\begin{equation}\label{eq5.18}
\begin{cases}
dX^{u^*}(t) = \frac{A(t)\mu(t)Z_0(t)}{Z(t)}dt + Z_0(t) dB(t); ; \; 0 \leq t \leq T\\
X^{u^*}(T)=E[(U')^{-1}(c \exp(-\int_0^T\theta(s)dB(s) -\frac{1}{2}\int_0^T (\theta(s))^2 ds))\Gamma(T)].
\end{cases}
\end{equation}

We have proved:

\begin{theorem}{(\emph{Optimal portfolio in an insider influenced market})}\\
The optimal portfolio $u^*$ for the problem \eqref{eq5.3} is given by 
\begin{equation}\label{eq5.17}
u^*(t) = \frac{Z_0(t)}{Z(t)},
\end{equation}
where $Z_0(t),Z(t)$ are the solutions of the BSDEs \eqref{eq5.1},\eqref{eq5.18}, respectively, and $c$ and $\theta$ are given by \eqref{eq5.16} and \eqref{eq5.9}, respectively.
\end{theorem}

\subsection{Predictive recursive utility maximization}
Consider a cash flow $X(t)=X^{c}(t)$ given by
\begin{equation}\label{eq5.20}
\begin{cases}
dX(t)= X(t)[\mu(t) dt + \sigma(t)dB(t)+\int_{R}\gamma(t,\zeta)\tN(dt,d\zeta)]-c(t) X(t) dt; \; t \in [0,T)\\
X(0) = x > 0.
\end{cases}
\end{equation}
Here $\mu(t),\sigma(t),\gamma(t,\zeta)$ are given bounded adapted processes, while $u(t):=c(t)$ is our control, interpreted as our \emph{relative consumption rate} from the cash flow. We say that $c$ is \emph{admissible} if $c$ is $\mathbb{F}$-adapted, $c(t) > 0$ and $X^c(t) > 0$ for all $t \in [0,T)$.
We put $\mathbb{G} = \mathbb{F}$.

Let $Y(t)=Y^c(t),Z(t)=Z^c(t),K(t,\zeta)=K^c(t,\zeta)$ be the solution of the predictive mean-field BSDE defined by
\begin{equation}\label{eq5.21}
\begin{cases}
dY(t)=- \{\alpha(t)A(t) + \ln (c(t)X(t))\} dt + Z(t) dB(t)+\int_{R}K(t,\zeta)\tN(dt,d\zeta); \; t \in [0,T)\\
Y(T) = 0,
\end{cases}
\end{equation}
where $\alpha(t) > 0$ is a given bounded $\mathbb{F}$-adapted process.
Then, inspired by classical definition of recursive utility in \cite{DE}, we define $Y^c(0)$ to be the \emph{predictive recursive utility} of the relative consumption rate $c$.

We now study the following predictive recursive utility maximization problem:

\begin{problem}
Find $c^* \in \mathcal{A}$ such that
\begin{equation}\label{eq5.22}
\sup_{c \in \mathcal{A}} Y^c(0) = Y^{c*}(0).
\end{equation}
\end{problem}

We apply the maximum principle to study this problem. In this case we have
$f=\varphi=h=0, \psi(x)=x$, and the Hamiltonian becomes
\begin{align}\label{eq5.23}
H(t,x,y,a,z,k,u,p,q,r,\lambda) &= x [(\mu(t) -c)p+\sigma(t) q+\int_{\mathbb{R}}\gamma(t,\zeta)r(\zeta)\nu(dt,d\zeta)] 
\nonumber \\ & \qquad +[a\alpha(t)+\ln c+\ln x]\lambda.
\end{align}
The associated backward-forward system of equations in the adjoint processes $p(t),q(t),\lambda(t)$ becomes
 \begin{itemize}
 \item{BSDE in $p(t),q(t)$:}
\begin{equation}\label{eq5.24}
\begin{cases}
dp(t) & = -[(\mu(t)-c(t))p(t)+\sigma(t)q(t) +\int_{\mathbb{R}} \gamma(t,\zeta)\nu(dt,d\zeta)+\frac{\lambda(t)}{X(t)}]dt\\
&+ q(t)dB(t) +\int_{\mathbb{R}}r(t,\zeta) \tN(dt,d\zeta) \; ; \; 0 \leq t \leq T \\
p(T) & = 0,
\end{cases}
\end{equation}
\item{SDE in $\lambda(t)$:}
\begin{equation}\label{eq5.25}
\begin{cases}
d\lambda(t) =  \alpha (t-\delta)\lambda(t-\delta)] \chi_{[\delta,T]}(t) dt\; ; \; 0 \leq t \leq T \\
\lambda(0) = 1.
\end{cases}
\end{equation}
\end{itemize}  
The delay SDE \eqref{eq5.25} does not contain any unknown parameters, and it is easily seen that it has a unique continuous solution $\lambda(t)>1$, which we may consider known.\\

We can now proceed along the same lines as in Section 5.2 of \cite{AO}: 
Maximizing $H$ with respect to $c$ gives the first order condition
\begin{equation}\label{eq5.26}
c(t) = \frac{\lambda(t)}{X(t)p(t)}.
\end{equation}  
The solution of the linear BSDE \eqref{eq5.24} is given by
\begin{equation}\label{eq5.27}
\Gamma(t)p(t)= E[\int_t^T \frac{\lambda(s)\Gamma(s)}{X(s)}ds | \mathcal{F}_t], 
\end{equation}
where
\begin{equation}
\begin{cases}
d\Gamma(t) = \Gamma(t^-)[(\mu(t)-c(t))dt +\sigma(t) dB(t) + \int_{\mathbb{R}} \gamma(t,\zeta)\tN(dt,d\zeta)] \; ; \; 0 \leq t \leq T \\
\Gamma(0)=1.
\end{cases}
\end{equation}
Comparing with \eqref{eq5.20} we see that
\begin{equation}\label{eq5.29}
X(t) = x \Gamma(t) \; ; \; 0 \leq t \leq T. \\
\end{equation}
Substituting this into \eqref{eq5.27} we obtain
\begin{equation}
p(t)X(t)= E[\int_t^T \lambda(s) ds | \mathcal{F}_t] \; ; \; 0 \leq t \leq T. \\
\end{equation}
Substituting this into \eqref{eq5.26} we get the following conclusion:

\begin{theorem}
The optimal relative consumption rate $c^*(t)$ for the predictive recursive utility consumption problem \eqref{eq5.22} is given by
\begin{equation}
c^*(t) = \frac{\lambda(t)}{E[\int_t^T \lambda(s) ds | \mathcal{F}_t]}\; ; \; 0 \leq t < T,
\end{equation}
where $\lambda(t)$ is the solution of the delay SDE \eqref{eq5.25}.
\end{theorem}

\end{document}